\documentclass[12pt]{amsart}
\usepackage{amscd,amssymb,amsmath,multicol}
\newtheorem{thm}[equation]{Theorem}
\numberwithin{equation}{section}
\newtheorem{cor}[equation]{Corollary}

\newtheorem{prop}[equation]{Proposition}
\newtheorem{tab}[equation]{Table}

\begin{document}
\raggedbottom \voffset=-.7truein \hoffset=0truein \vsize=8truein
\hsize=6truein \textheight=8truein \textwidth=6truein
\baselineskip=18truept
\def\cal{\mathcal}
\def\vareps{\varepsilon}
\def\mapright#1{\ \smash{\mathop{\longrightarrow}\limits^{#1}}\ }
\def\mapleft#1{\smash{\mathop{\longleftarrow}\limits^{#1}}}
\def\mapup#1{\Big\uparrow\rlap{$\vcenter {\hbox {$#1$}}$}}
\def\mapdown#1{\Big\downarrow\rlap{$\vcenter {\hbox {$\ssize{#1}$}}$}}
\def\on{\operatorname}
\def\span{\on{span}}
\def\cat{\on{cat}}
\def\a{\alpha}
\def\bz{{\Bbb Z}}
\def\gd{\on{gd}}
\def\imm{\on{imm}}
\def\sq{\on{Sq}}
\def\sspan{\on{span}^0}
\def\rank{\on{rank}}
\def\eps{\epsilon}
\def\br{{\Bbb R}}
\def\bc{{\Bbb C}}
\def\bh{{\Bbb H}}
\def\tfrac{\textstyle\frac}
\def\w{\wedge}
\def\b{\beta}
\def\A{{\cal A}}
\def\P{{\cal P}}
\def\zt{{\Bbb Z}_2}
\def\bq{{\Bbb Q}}
\def\ker{\on{ker}}
\def\coker{\on{coker}}
\def\u{{\cal U}}
\def\e{{\cal E}}
\def\exp{\on{exp}}
\def\wbar{{\overline w}}
\def\xbar{{\overline x}}
\def\ybar{{\overline y}}
\def\zbar{{\overline z}}
\def\ebar{{\overline e}}
\def\nbar{{\overline n}}
\def\mbar{{\overline m}}
\def\ubar{{\overline u}}
\def\et{{\widetilde E}}
\def\ni{\noindent}
\def\coef{\on{coef}}
\def\den{\on{den}}
\def\gd{{\on{gd}}}
\def\ev{{\text{ev}}}
\def\od{{\text{od}}}
\def\N{{\Bbb N}}
\def\Z{{\Bbb Z}}
\def\Q{{\Bbb Q}}
\def\R{{\Bbb R}}
\def\C{{\Bbb C}}
\def\Bin{\on{Bin}}
\title[Vector fields]
{Vector fields on $RP^m\times RP^n$}

\author{Donald M. Davis}
\address{Department of Mathematics, Lehigh University\\Bethlehem, PA 18015, USA}
\email{dmd1@lehigh.edu}
\date{December 15, 2010}


\maketitle
\begin{abstract} The span of a manifold is its maximum number of linearly independent vector fields.
We discuss the question, still unresolved, of whether $\span(P^m\times P^n)$ always equals $\span(P^m)+\span(P^n)$.
Here $P^n$ denotes real projective space.
We use $BP$-cohomology to obtain new upper bounds for $\span(P^m\times P^n)$, much stronger than previously known bounds.
 \end{abstract}
\section{Introduction}\label{intro}
The {\bf span} of a manifold is the maximum number of linearly independent vector fields on it. The following result
is well-known. Throughout the paper, $\nu(n)$ denotes the exponent of 2 in $n$, and $P^n$ denotes real projective space.
\begin{thm}  Let $V(n)=8a+2^b-1$ if $\nu(n+1)=4a+b$ with $0\le b\le 3$. Then $\span(P^n)=V(n)$.
\end{thm}
\begin{proof} It was proved by Adams in \cite{vfld} that $\span(S^n)=V(n)$. Since linearly independent (l.i.)
vector fields on $P^n$ pull back to l.i.~vector fields on $S^n$, this implies $\span(P^n)\le V(n)$.
Since the vector fields on $S^n$ can be chosen to satisfy $v_i(-x)=-v_i(x)$, they induce l.i.~vector fields on $P^n$. (\cite[p.140]{Huse})
\end{proof}
This of course implies that
\begin{equation}\span(P^m\times P^n)\ge V(m)+V(n).\label{vv}
\end{equation}
Although it seems unlikely that equality always holds in (\ref{vv}), there are no known examples
in which $\span(P^m\times P^n)$ exceeds $V(m)+V(n)$. Our first result
shows that equality does hold in (\ref{vv}) in many cases.
\begin{thm}\label{equ} If neither $m+1$ nor $n+1$ is divisible by 16, or if $m=1$, $3$, or $7$, then
$$\span(P^m\times P^n)=V(m)+V(n).$$
\end{thm}
\noindent The proof of this, which is quite elementary, is given in Section \ref{pfs}.

Our second result uses $BP$-cohomology to obtain new upper bounds for $\span(P^m\times P^n)$
which are exponentially stronger than previously known results. This should be considered the main
result of the paper.
\begin{thm}\label{bnd} Let $r=\nu(M)\ge4$ and $t=\nu(N)\le r$. Then \begin{eqnarray*}&&\span(P^{2M-1}\times P^{2N-1})
\le\begin{cases}14\cdot2^e-4&\\
18\cdot2^e-4&\\
20\cdot2^e+8k+6\\
24\cdot2^e+4k-2&\end{cases}\\
&&\text{if}\quad r+t=2^e+2e+\begin{cases}2&\\ 3&\\ 4+k,\quad 0\le k<2^e-2\\ 4+k,\quad k=2^e-2,2^e-1&\end{cases}\text{and}\quad t>
\begin{cases}e&\\e+1&\\e&\\e.\end{cases}\end{eqnarray*}
Here $e\ge1$ or $(r,t,e,k)=(4,1,0,0)$.
\end{thm}
We prove this result in Section \ref{pfs}. In Section \ref{nums}, we give numerical illustrations of this theorem
and compare it with previous results. We also discuss some cases in which it can be extended.

\section{Proofs} \label{pfs}
In this section, we prove Theorems \ref{equ} and \ref{bnd}.
The following result is well-known.
\begin{prop} \label{SWprop}$\span(P^m\times P^n)\le 2^{\nu(m+1)}+2^{\nu(n+1)}-2$.\end{prop}
\begin{proof} The total Stiefel-Whitney class of $P^m\times P^n$ is $(1+x_1)^{m+1}(1+x_2)^{n+1}$, truncated after
$x_1^m$ and $x_2^n$. By well-known properties of binomial coefficients mod 2, the highest nonzero Stiefel-Whitney
class is $x_1^{m+1-2^{\nu(m+1)}}x_2^{n+1-2^{\nu(n+1)}}$. Thus the tangent bundle cannot be stably equivalent
to a bundle of dimension less than $$d:=m+1-2^{\nu(m+1)}+n+1-2^{\nu(n+1)}.$$ Hence the number of l.i.~vector fields is at most
$m+n-d$,
as claimed.\end{proof}

Now we can prove our first theorem.
\begin{proof}[Proof of Theorem \ref{equ}.]
The first case follows from the previous proposition together with the fact that $V(n)=2^{\nu(n+1)}-1$
if $n\not\equiv15$ mod 16. For the second case, we need the important notion of {\it stable span}.
The span of a vector bundle $\theta$ is its maximal number of l.i.~sections, and
the stable span of a manifold $M$, denoted $\sspan(M)$, equals $\span(\tau(M)+m\eps)-m$ for $m>0$, which is easily seen
to be independent of such $m$. Here $\tau$ denotes the tangent bundle, and $\eps$ a trivial bundle.

The restriction of $\tau(P^m\times P^n)$ to $*\times P^n$ is $\tau(P^n)+m\eps$. Hence $$\span(P^m\times P^n)\le m+\sspan(P^n).$$
For $m=1$, 3, or 7, $m=\span(P^m)$ and so it remains to show that $\sspan(P^n)=\span(P^n)$.

If $n$ is even, both are 0 since $w_n(\tau(P^n))\ne0$. In \cite[1.11]{JT}, James and Thomas proved that if $n$ is odd, an $n$-plane bundle over
$P^n$ is equivalent to $\tau(P^n)$ if and only if they are stably equivalent, implying the result in this case. \end{proof}

In the rest of this section, we prove our second theorem, which is more substantial.
 Since $\tau(P^{2M-1}\times P^{2N-1})+2\eps\approx 2M\xi_{2M-1}\times 2N\xi_{2N-1}$, then,
using the $BP$-Euler class as in \cite[pp.331-332]{SW}, we obtain that if $\span(P^{2M-1}\times P^{2N-1})>s$, then
\begin{equation}\label{bc}\sum_{i,j}(-1)^{i+j}\tbinom Mi\tbinom Nj x_1^i x_2^j x_3^{M+N-i-j}=0\in BP^{2M+2N}(P^{2M-2}\times P^{2N-2}\times P^{s+2}).\end{equation}
Note that we have restricted to $P^{2M-2}\times P^{2N-2}\subset P^{2M-1}\times P^{2N-1}$ to simplify the calculation.
Here $x_i$ is a 2-dimensional class corresponding to the $i$th factor, and $BP$ is the 2-local Brown-Peterson spectrum. 
The conclusion (\ref{bc}) also holds with
Cartesian product replaced by smash product, using the direct sum splitting.

It will be convenient to work with the Johnson-Wilson spectrum $BP\langle3\rangle$. See, for example, \cite[p.117]{JWY}.
We will call if $B$.
Recall that $B_*=\bz_{(2)}[v_1,v_2,v_3]$ with $|v_i|=2^{i+1}-2$. Let $v_0=2$ and let $I$ denote the ideal
$(v_0,v_1,v_2,v_3)$. There is a power series $[2](x)=\sum a_jx^{j+1}$ with $a_j\in B_{2j}$
satisfying $a_0=2$ and $$[2](x)\equiv v_0x+v_1x^2+v_2x^4+v_3x^8 \mod I^2.$$
See, for example, \cite[p.120]{JWY}.
Let $P_1=P^\infty$ and $B^*=B_{-*}$. Gysin sequence arguments show that  $$B^*(P^{2n-2})\approx B^*[x]/(x^{n},[2](x))$$ and $$B_*(P_1)\approx B_*(z_i:i\ge1)/(\sum a_jz_{i-j}),$$
where $|z_i|=2i-1$, and there are duality isomorphisms
\begin{equation}\label{Piso}B^{2n-2i}(P^{2n-2})\approx B_{2i-1}(P_1)\end{equation}
 for $i<n$ under which
$x^{n-i}$ corresponds to $z_i$. These isomorphisms pass to 2- and 3-fold smash products, as described briefly in
\cite[p.330]{SW}. If $i$, $j$, and $k$ are positive integers, let $[i,j,k]\in B_{2(i+j+k)-3}(P_1\w P_1\w P_1)$
denote the external product of classes $z_i$, $z_j$, and $z_k$, as in \cite[p.120]{JWY}.
Let $$Q=B_*(P_1)\otimes_{B_*}B_*(P_1)\otimes_{B_*}B_*(P_1).$$
By \cite[1.3]{JWY}, $Q$ is a subgroup of $B_*(P_1\w P_1\w P_1)$. This $Q$ is a $B_*$-module. Let $F_s=I^s\cdot Q\subset Q$. This filtration could be thought of as filtration in an Adams spectral sequence. The content of
\cite[2.3]{JWY} can be restated as follows.
\begin{prop}\label{JWYprop} For all $s\ge0$, $F_s/F_{s+1}$ is a graded $\zt$-vector space with basis all $v_3^s[i_1,i_2,i_3]$, $i_j>0$.\end{prop}

The action of $v_0:F_0/F_1\to F_1/F_2$ can be determined using the 2-series $[2](x)$ in the following way.
Temporarily write $[e_1,e_2,e_3]$ as $z_1^{e_1}z_2^{e_2}z_3^{e_3}$. The 2-series forces relations
$$(v_0+v_1z_i^{-1}+v_2z_i^{-3}+v_3z_i^{-7})z_1^{e_1}z_2^{e_2}z_3^{e_3}=0\text{ in }F_1/F_2.$$
We apply the relation to repeatedly replace $v_0$ by $v_1z_1^{-1}+v_2z_1^{-3}+v_3z_1^{-7}$, $v_1$ by $v_0z_2+v_2z_2^{-2}+v_3z_2^{-6}$,
and $v_2$ by $v_0z_3^3+v_1z_3^2+v_3z_3^{-4}$. This is the approach that was taken in \cite{SW}. For $i=0$, 1, and 2, we obtain infinite series
$f_i(z_1,z_2,z_3)$ such that $v_i=v_3\cdot f_i(z_1,z_2,z_3)$. For example, the determination of $f_0$ begins as follows.
\begin{eqnarray*} v_0&=&v_1z_1^{-1}+v_2z_1^{-3}+v_3z_1^{-7}\\
&=&(v_0z_2+v_2z_2^{-2}+v_3z_2^{-6})z_1^{-1}\\
&&+(v_0z_3^3+v_1z_3^2+v_3z_3^{-4})z_1^{-3}+v_3z_1^{-7}\\
&=&(v_1z_1^{-1}+v_2z_1^{-3}+v_3z_1^{-7})(z_1^{-1}z_2+z_1^{-3}z_3^3)\\
&&+(v_0z_2+v_2z_2^{-2}+v_3z_2^{-6})z_1^{-3}z_3^2\\
&&+(v_0z_3^3+v_1z_3^2+v_3z_3^{-4})z_1^{-1}z_2^{-2}\\
&&+v_3(z_1^{-7}+z_1^{-1}z_2^{-6}+z_1^{-3}z_3^{-4})\\
&=&\cdots\end{eqnarray*}
The three monomials in the last line above could be thought of as the start of $f_0$. However, they may (and will)  be cancelled later
in the algorithm. The procedure does converge in the sense that any monomial can only appear a finite number of times. This is true
because every step decreases $3\nu_{z_1}+\nu_{z_2}-2\nu_{v_0}+\nu_{v_2}$ in all monomials.

The series $f_0$, $f_1$, and $f_2$ must satisfy the equations
\begin{eqnarray*}f_0&=&f_1z_1^{-1}+f_2z_1^{-3}+z_1^{-7}\\
f_1&=&f_0z_2+f_2z_2^{-2}+z_2^{-6}\\
f_2&=&f_0z_3^3+f_1z_3^2+z_3^{-4}.\end{eqnarray*}
Rearrange them so that they appear as a system over $\zt$ for the three unknowns $f_0$, $f_1$, and $f_2$, and solve by Cramer's
rule to obtain
\begin{eqnarray*}f_0&=&\frac{z_1^{-1}z_2^{-2}z_3^{-4}+z_1^{-1}z_2^{-6}+z_1^{-3}z_2^{-6}z_3^2+z_1^{-3}z_3^{-4}+z_1^{-7}+z_1^{-7}z_2^{-2}z_3^2}
{1+z_2^{-2}z_3^2+z_1^{-1}z_2+z_1^{-1}z_2^{-2}z_3^3+z_1^{-3}z_2z_3^2+z_1^{-3}z_3^3}\\
&=&z_1^{-1}z_2^{-2}z_3^{-4}+z_1^{-1}z_2^{-4}z_3^{-2}+z_1^{-2}z_2^{-1}z_3^{-4}+z_1^{-2}z_2^{-4}z_3^{-1}+z_1^{-4}z_2^{-1}z_3^{-2}+z_1^{-4}z_2^{-2}z_3^{-1}.
\end{eqnarray*}
The latter equation can be verified by cross multiplication. Note that it turns out that $f_0$ is not an infinite series after all.
Our conclusion is stated in the next two results.
\begin{thm} For distinct nonnegative integers $i$, $j$, and $k$, let
$\cal{P}(i,j,k)$ denote the set consisting of the six permutations of $i$, $j$, and $k$.
Then, in $Q$,
$$2[e_1,e_2,e_3]\equiv v_3\sum_{(i,j,k)\in{\cal P}(1,2,4)}[e_1-i,e_2-j,e_3-k]\mod F_2.$$
\end{thm}
Iterating this, we obtain
\begin{cor}\label{v3cor} First,
$$2^{2^t}[e_1,e_2,e_3]\equiv v_3^{2^t}\sum[e_1-i,e_2-j,e_3-k]\mod F_{2^t+1},$$
where the sum is taken over all $(i,j,k)\in{\cal P}(2^t,2^{t+1},2^{t+2})$. More generally, if $m=\sum\limits_{\ell=1}^d 2^{t_\ell}$ with $\{t_\ell\}$
distinct, then
$$2^m[e_1,e_2,e_3]\equiv v_3^m\sum[e_1-i_1-\cdots-i_d,e_2-j_1-\cdots-\j_d,e_3-k_1-\cdots-k_d]\mod F_{m+1},$$
summed over all $(i_\ell,j_\ell,k_\ell)\in{\cal P}(2^{t_\ell},2^{t_\ell+1},2^{t_\ell+2})$ with $\ell=1,\ldots,d$.
\end{cor}

 The implication of $\span(P^{2M-1}\times P^{2N-1})>2s'-4$
in (\ref{bc}) becomes
$$\sum(-1)^{i+j}\tbinom M{M-i}\tbinom N{N-j}x_1^{M-i}x_2^{N-j}x_3^{s'-(s'-i-j)}=0\in B^*(P^{2M-2}\w P^{2N-2}\w P^{2s'-2})$$
and then, under the isomorphism (\ref{Piso}),
\begin{equation}\sum(-1)^{i+j}\tbinom Mi\tbinom Nj[i,j,s'-i-j]=0\in Q.\label{ijk}\end{equation}

We preview the detailed proof of Theorem \ref{bnd} by illustrating with the case $r=5$, $t=3$.
This is the fourth case in the theorem, with $e=1$ and $k=0$.
The claim then is that if
$M\equiv32$ mod 64 and $N\equiv8$ mod 16, then $\span(P^{2M-1}\times P^{2N-1})\le46$. To prove this, we assume
$\span(P^{2M-1}\times P^{2N-1})>46$ and deduce as above that \begin{equation}\label{S25}\sum(-1)^{i+j}\tbinom{64k+32}i\tbinom{16\ell+8}j[i,j,25-i-j]=0\in Q.\end{equation}
Since all nonzero terms have $i+j\le 24$, all terms in the sum are divisible by $2^2$, and so our sum is in $F_2$.
Here and later we use that $\nu\binom {u2^r}i\ge r-\nu(i)$,
with equality if $i\le 2^r$.
Note that the filtration-1 term $2[16,8,1]$ which occurs in (\ref{S25}) is 0 since subtracting a permutation of $(4,2,1)$ from $(16,8,1)$
always results in a non-positive entry.

The only terms in (\ref{S25}) not divisible by $2^3$ are
 $2^2u_1[16,4,5]$ and $2^2u_2[8,8,9]$ with $u_i$ odd, and so, mod $F_3$, our sum equals $2^2[16,4,5]+2^2[8,8,9]$.
By Corollary \ref{v3cor}, this is equal, mod $F_3$, to
$$v_3^2([16-8,4-2,5-4]+[8-2,8-4,9-8]+[8-4,8-2,9-8])=v_3^2([8,2,1]+[6,4,1]+[4,6,1]),$$
which is nonzero by \ref{JWYprop}, contradicting the assumption that the span is greater than 46.

Now we begin the proof of Theorem \ref{bnd} in earnest. For the first case in \ref{bnd}, we assume $\span(P^{2M-1}\times P^{2N-1})>14\cdot2^e-4$. Then, with $u_i$ odd,
\begin{equation}\label{sum1}\sum_{i,j}\pm\binom{u_12^r}i\binom{u_22^t}j[i,j,7\cdot2^e-i-j]=0.\end{equation}
 Here the pair $(r,t)$ appears in the list
$$(2^e+e+1,e+1),\ldots,(2^{e-1}+e+1,2^{e-1}+e+1),$$
always summing to $2^e+2e+2$. Since we must have $i+j<7\cdot 2^e$, the terms in the sum having lowest filtration are
$T_1:=2^{2^e-1}[2^{e+2},2^{e+1},2^e]$ and, if $t>e+1$, $T_2:=2^{2^e-1}[2^{e+1},2^{e+2},2^e]$.
By Corollary \ref{v3cor}, the
term $T_1$ equals, mod $F_{2^e}$,
 $$v_3^{2^e-1}[2^{e+2}-(2^{e+2}-4),2^{e+1}-(2^{e+1}-2),2^e-(2^e-1)]=v_3^{2^e-1}[4,2,1]$$
plus perhaps $v_3^{2^e-1}$ times other terms of degree $2\cdot7$. The  term $T_2$ equals, mod $F_{2^e}$, the sum of terms $$v_3^{2^e-1}[2^{e+1}-\sum_{i=0}^{e-1} 2^{i+a_i},2^{e+2}-\sum_{i=0}^{e-1} 2^{i+b_i},2^e-\sum_{i=0}^{e-1} 2^{i+c_i}],$$
where each $(a_i,b_i,c_i)$ is a permutation of $(0,1,2)$.
We must have $\sum 2^{i+b_i}\le 2^{e+2}-4$, and hence $T_2$ has no terms of the form $v_3^{2^e-1}[-,2,-]$.
Thus the nonzero term $v_3^{2^e-1}[4,2,1]$ in $T_1$ is uncancelled in filtration $2^e-1$, and so the LHS of (\ref{sum1})
is nonzero, a contradiction.

For the second case, if span $>18\cdot2^e-4$, then
\begin{equation}\label{sum2}\sum\pm\binom{u_12^r}i\binom{u_22^t}j[i,j,9\cdot2^e-i-j]=0.\end{equation}
Under the hypotheses, the only term of lowest filtration is $2^{2^e-1}[2^{e+2},2^{e+2},2^e]$, and this
equals, mod $F_{2^e}$, $v_3^{2^e-1}[2^{e+1}+2,4,1]$ plus  other terms, and hence is nonzero.

The third and fourth cases are distinguished by which has nonzero terms in the smaller grading.
For the fourth case, the terms of smallest filtration are $$2^{2^e+k}[2^{e+3},2^{e+1},2^{e+1}+1+2k],$$ possibly $2^{2^e+k}[2^{e+1},2^{e+3},2^{e+1}+1+2k]$, and if $k=0$ then also
 $2^{2^e}[2^{e+2},2^{e+2},2^{e+2}+1]$.
The first term equals, mod $F_{2^e+k+1}$, $v_3^{2^e+k}[2^{e+2}-4k,2^e-k,1]$ plus possibly other terms, and this cannot be cancelled by either
of the others.

Under the hypothesis $k<2^e-2$ of the third case, there is a term $$2^{2^e+k+1}[2^{e+2},2^{e+1},2^{e+2}+4k+5]$$ of smaller grading than  the terms just considered, and in this
grading there are no terms of smaller filtration.
If $t>e+1$, there is also a nonzero term $$2^{2^e+k+1}[2^{e+1},2^{e+2},2^{e+2}+4k+5].$$ The first term equals, mod higher filtration and
other terms, $$T:=v_3^{2^e+k+1}[2^{e+1}-2k-2,2^e-k-1,1].$$ When the second term is rewritten as a sum of terms $$v_3^{2^e+k+1}[2^{e+1}-\sum i_\ell,2^{e+2}-\sum j_\ell,2^{e+2}+4k+5-\sum k_\ell]$$ as in \ref{v3cor}, the only way it could contain the term  $T$ is if $\sum k_\ell=
2^{e+2}+4k+4$. Then $\sum j_\ell\le 2^{e+1}+2k+2$, and since $k<2^e-2$, we must have $2^{e+2}-\sum j_\ell>2^e-k-1$, and so $T$ cannot be cancelled.

\section{Numerical results}\label{nums}
Theorem \ref{bnd} is exponentially better than the previous best known results, most of which are given in Proposition
\ref{SWprop}. For a typical example, if $\nu(m+1)=\nu(n+1)=2^{e-1}+e+2$, then we obtain
$$\span(P^m\times P^n)\le 2^{e+4}-2^{e+1}-4,$$
while the result from \ref{SWprop} is
$$\span(P^m\times P^n)\le 2^{2^{e-1}+e+3}-2.$$

In \cite{Suz}, $K$-theoretic methods were used to obtain nonexistence results for vector fields on
products of real projective spaces. In \cite{Kob}, slight improvements were obtained in some cases, but not in any applicable to this paper.
All these results are weaker than those of Proposition \ref{SWprop} unless
both $m+1$ and $n+1$ are 2-powers, or $m+1$ is a 2-power and $n<m/2$. The $K$-theory bound is always more than
$1/2$ times the Stiefel-Whitney bound. See Table \ref{tab1} for numerical examples. In our example in the previous paragraph, the $K$-theoretic
methods give no new information, compared to Proposition \ref{SWprop}, except in the case of $P^m\times P^m$ with $m=2^{2^{e-1}+e+2}-1$, and in this case the bound is roughly $3\cdot2^{2^{e-1}+e+1}+2^{e-1}+e-1$,
still exponentially larger than our bound.

Theorem \ref{bnd} can be extended to include cases in which $t$-values are smaller than they are in that theorem.
We did not include them there because their patterns become too complicated. We won't even list many of them here because,
for best results, they begin to involve more than just the 2-divisibility of $N$.
The following result describes what we can deduce when $t$ is 1 (resp. 1 or 2) smaller than is allowed in the first (resp. second)
case of Theorem \ref{bnd}.
\begin{prop} \label{other}Let $e\ge1$, $r=\nu(M)$, and $t=\nu(N)$. Then
\begin{eqnarray*}&&\span(P^{2M-1}\times P^{2N-1})\\
&\le&\begin{cases}20\cdot2^e-2&\text{if }r=2^e+e+2\text{ and }t=e+1\text{ or }(e>1\text{ and }N\equiv3\cdot2^e\ {\rm mod}\ 2^{e+2})\\ 22\cdot2^e-6&\text{if }r=2^e+e+2\text{ and }t=e=1\text{ or }N\equiv2^e\ {\rm mod}\ 2^{e+2}\\ 20\cdot2^e+6&\text{if }r=2^e+e+3\text{ and }N\equiv3\cdot2^e\ {\rm mod}\ 2^{e+2}\\ 38\cdot2^e-6&\text{if }r=2^e+e+3\text{ and }t=e.\end{cases}
\end{eqnarray*}
\end{prop}
Note that the third case here is a strengthening of the fourth case which applies when more information about
$N$ is known other than just its 2-exponent. Other similar strengthenings can be given, but become too tedious to list.
\begin{proof} These are proved by the same method as the proof of \ref{bnd}. The determining terms of lowest filtration in the
four cases here are
\begin{eqnarray*}2^{2^e}[2^{e+2},2^{e+1},2^{e+2}+1]&\sim &v_3^{2^e}[2^{e+1},2^e,1]\\
2^{2^e-1}[2^{e+3},2^e,2^{e+1}-1]&\sim&v_3^{2^e-1}[2^{e+2}+4,1,1]\\
2^{2^e+1}[2^{e+2},2^{e+1},2^{e+2}+5]&\sim&v_3^{2^e+1}[2^{e+1}-2,2^e-1,1]\\
2^{2^e-1}[2^{e+4},2^e,2^{e+1}-1]&\sim&v_3^{2^e-1}[3\cdot2^{e+2}+4,1,1].\end{eqnarray*}
\end{proof}

We close with a table comparing the various bounds for $\span(P^{m}\times P^{111})$ for $m=2^e-1$ and $m=3\cdot2^e-1$
for $5\le e\le17$. The first column (after the $e$-column) gives the lower bound for $\span(P^m\times P^{111})$ given
by (\ref{vv}). It is the same for $m=2^e-1$ and $m=3\cdot2^e-1$. The second column gives our upper bound for $\span(P^m\times P^{111})$ from \ref{bnd}, \ref{other}, and a slight extension of the third case of \ref{other}. In this range, it is the
same for either value of $m$.
The third column gives the upper bound given by Stiefel-Whitney classes (Proposition \ref{SWprop}). It is always the same
for the two $m$-values. The final column gives the upper bound obtained by Suzuki in \cite{Suz}  when $m=2^e-1$.
If $m=3\cdot2^e-1$, his bound is larger than the Stiefel-Whitney bound, and so we do not bother to list it.

\medskip
\begin{minipage}{6.5in}
\begin{tab}\label{tab1}{Bounds for $\span(P^m\times P^{111})$, $m=2^e-1$ or $3\cdot2^e-1$}
\begin{center}
\begin{tabular}{r|cccc}
&&&&Suzuki\\
&&our&Stief-Whit&upper\\
&lower&upper&upper&bound\\
$e$&bound&bound&bound&$m=2^e-1$\\
\hline
$5$&$17$&$32$&$46$&$130$\\
$6$&$19$&$46$&$78$&$148$\\
$7$&$23$&$50$&$142$&$182$\\
$8$&$24$&$52$&$270$&$246$\\
$9$&$25$&$78$&$526$&$374$\\
$10$&$27$&$86$&$1038$&$630$\\
$11$&$31$&$94$&$2062$&$1146$\\
$12$&$32$&$102$&$4110$&$2170$\\
$13$&$33$&$106$&$8206$&$4218$\\
$14$&$35$&$158$&$16398$&$8316$\\
$15$&$39$&$166$&$32782$&$16510$\\
$16$&$40$&$174$&$65550$&$32894$\\
$17$&$41$&$182$&$131086$&$65662$
\end{tabular}
\end{center}
\end{tab}
\end{minipage}
\medskip

Note that our upper bound is moderately close to the known lower bound. For fairly large $e$, on the other hand, the other known upper bounds are
exponentially larger than ours and the lower bound.

\def\line{\rule{.6in}{.6pt}}


\begin{thebibliography}{99}

\bibitem{vfld} J.F.Adams, {\em Vector fields on spheres}, Annals of Math {\bf 75} (1962) 603-632.
\bibitem{Huse} D.Husemoller, {\it Fiber bundles}, Springer-Verlag (1993).
\bibitem{JT} I.M.James and E.Thomas, {\em An approach to the enumeration problem for non-stable vector bundles}, Jour Math Mech {\bf 14} (1965) 485-506.
\bibitem{JWY} D.C.Johnson, W.S.Wilson, and D.Y.Yan, {\em Brown-Peterson homology of elementary abelian $p$-groups, II}, Topology and its Appls {\bf 59} (1994) 117-136.
\bibitem{Kob} T.Kobayashi, {\em Note of $\gamma$-dimension and products of real projective spaces}, Jour Math Soc Japan {\bf 34} (1982) 501-505.
\bibitem{SW} H.-J.Song and W.S.Wilson, {\em On the nonimmersion of products of real projective spaces}. Trans Amer Math Soc {\bf 318} (1990) 327-334.
\bibitem{Suz} H.Suzuki, {\em Operations in $KO$-theory and products of real projective spaces}, Mem Fac Sci Kyushu Univ {\bf 18} (1964) 140-153.
\end{thebibliography}
\end{document}